\documentstyle[11pt,twoside]{article}

\def\R{\relax\ifmmode I\!\!R\else$I\!\!R$\fi}

\def\Z{\relax\ifmmode Z\!\!\!Z\else$Z\!\!\!Z$\fi}

\def\C{\relax\ifmmode C\!\!\!\!I\else$C\!\!\!\!I$\fi}

\def\K{\relax\ifmmode I\!\!K\else$I\!\!K$\fi}

\def\N{\relax\ifmmode I\!\!N\else$I\!\!N$\fi}

\def\colon{{:}\;}
\newcounter{defcounter}[section]
\newenvironment{remark}%
{\vspace{0.2cm}\begin{sloppypar}\noindent\stepcounter{defcounter}{\bfseries Remark
      \thesection.\thedefcounter}}%
{\end{sloppypar}\vspace{0.2cm}} 
\newtheorem{theorem}{Theorem}[section]
\newtheorem{proposition}{Proposition}[section]

\newcommand{\proof}{{\noindent\bf Proof }}

\newcommand{\qed}{\hfill $\Box$}

\begin{document}
\parindent=1em
\baselineskip 15pt
\hsize=12.3 cm
\vsize=18.5 cm 

\thispagestyle{empty}
\begin{center}
{\Large {\bf Number theoretical peculiarities in the dimension theory 
of dynamical systems}}
\end{center}
\begin{center}
J. Neunh\"auserer\footnote{Supported by ''DFG-Schwerpunktprogramm - Dynamik: Analysis, effiziente Simulation und Ergodentheorie''.}\\
Fachbereich Mathematik, Technische Universit\"at Dresden\\
Zellescher Weg. 12-14, 01069 Dresden, Germany\\
e-mail: neuni@math.tu-dresden.de
\end{center}
\begin{abstract}We show that dimensional theoretical properties
of dynamical systems can considerably change because of number theoretical
peculiarities of some parameter values
\\{\bf AMS subject classification 2000}: 37C45, 37A45
\end{abstract}
\section{Introduction \label{Chap1}}
In the last decades there has been an enormous interest in geometrical invariants
of dynamical systems especially in the Hausdorff dimension of invariant sets
like attractors, reppelers or hyperbolic sets
and ergodic measures on these sets. A dimension theory of dynamical systems
was developed and now a days the Hausdorff dimension seems to have
its place beside classical invariants like entropy or Lyapunov exponents.\footnote{We refer appendix A of this work and to the book of Falconer \cite{[FA1]} for an introduction to dimension theory and recommend the book of Pesin \cite{[PE1]} for the dimension theory of dynamical systems.} \newline
There are two main principles that form a kind of a guide line
through the dimension theory of dynamical systems. The first 
states the identity of Hausdorff and box-counting dimension of invariant 
sets. The second one
is the variational principle for Hausdorff dimension which states that
the Hausdorff dimension of a given invariant set can be approximated
by the Hausdorff dimension of ergodic measures on these set or in a
stronger form states the existence of an ergodic measure of
full Hausdorff dimension on an given invariant set. In many situations
these principle are essential to determine the 
Hausdorff dimension of an invariant set and for relating this quantity to
other characteristics of the dynamics like entropy, 
Lyapunov exponents and pressure.\newline
For conformal repellers we know that the identity
of Hausdorff and box-counting dimension holds and that there exists an ergodic measure of full Hausdorff dimension (see chapter 7 of \cite{[PE1]}). For 
hyperbolic sets of diffeomorphisms the variational principle
for Hausdorff dimension does not hold in general (see \cite{[MM]}). But again 
if the system is conformal restricted to stable resp. unstable manifolds
there exits an ergodic measure of full dimension for the restrictions and the
identity of box-counting and Hausdorff dimension of the hyperbolic set holds (see again chapter 7 of \cite{[PE1]}).
In the non conformal situation there is no general theory this days that
allows us to determine the dimensional theoretical properties of a given dynamical system.
But there are a lot of results for special classes of systems that state
that the variational principle or the identity of box-counting and Hausdorff
dimension or both hold at least generically in the sense of Lebesgue measure
on the parameter space (see for instance \cite{[FA2]}, \cite{[SCH]}, \cite{[NE]}, \cite{[SS]}, \cite{[SO2]}). In this paper we focus at such classes of systems.\\
We will show that in situation were there generically
exists an ergodic measure of full Hausdorff dimension the variational principle for Hausdorff dimension may not hold in general because of
number theoretical peculiarities of some parameter values (see Theorem 2.1 below). Furthermore we will show that the identity of Box-Counting and Hausdorff dimension may drop because of number theoretical peculiarities in situations were this identity generically holds (see Theorem 2.2 below). Our example for the first phenomena is the Fat Baker's transformation
and our example for the second phenomena is a class of self-affine reppelers.
Both classes of systems are very simple but it seems obvious to us that the same phenomena appear as well in more complicated examples also this would be of course even harder to proof.   \newline
All our results are related to a special class of
algebraic integers namely Pisot-Vijayarghavan numbers\footnote{A  Pisot-Vijayarghavan number is an algebraic integer with all
its algebraic conjugates inside the unit circle (see  appendix B)} (short: PV numbers) 
and they are in some sense the consequence of a generalisation of results of Erd\"os \cite{[ER1]}, Garsia (\cite{[GA1]}, \cite{[GA2]}) and 
Alexander and Yorke \cite{[AY]} on the singularity and dimension of invenitly convolved measures. We think that from the 
viewpoint of geometric measure theory and algebraic number theory this generalisation is interesting in itself (see Theorem 4.1 below).\newline
The rest of the paper is organised as follows. 
In section two we define the systems we study, state our main Theorems 2.1 and 2.2 about these systems and comment on our results. In section three we introduce coding maps for our systems and find representations of all ergodic measures using these codings. In section four we
define a class of Borel probability measures associated with a PV numbers (Erd\"os measures), introduce a kind of entropy 
related to this measure (Garsia entropy) and state our main Theorem 4.1 about the singularity and the Hausdorff dimension
of Erd\"os measures.
The proof of Theorem 4.1 is given in section five, the proof of Theorem 2.1 is contained in section
six and the proof of Theorem 2.2 can be found in section seven. All our proofs consist of several propositions
which may be interesting in them self.
In appendix A we collect some basic definitions and facts in dimension theory 
and in appendix B we we define PV numbers and state the properties of these algebraic integers that we need in our work. \\\newline
{\bf Acknowlegdment:} I wish to thank J\"org Schmeling who helped me a 
lot to find the results presented here.  
\section{Basic definitions and main results}
For $\beta\in (0.5,1)$ we define the {\bf Fat Baker's transformation} 
$f_{\beta}\colon \R \times [-1,1] \to \R \times [-1,1]$ by  
\[  
{ f_{\beta} (x,y) } = \lbrace
\begin{array}{cc}
(\beta x+(1-\beta),2y-1)\quad \mbox{if}\quad y \geq 0   \\
(\beta x-(1-\beta),2y+1)\quad \mbox{if} \quad y<0.    
\end{array}
\]
This map was introduced by Alexander and Yorke in \cite{[AY]}. It is 
called Fat Baker's transformation because if we set $\beta=0.5$ we get the classical Baker's transformation.\newline
It is obvious that the attractor of $f_{\beta}$ is the whole square $[-1,1]^{2}$ which has Hausdorff and box-counting dimension two. We always restrict $f_{\beta}$ to its attractor.\newline\newline
Now we state our main result about the Fat Baker's transformation.
\begin{theorem}
If $\beta\in (0.5,1)$ is the reciprocal of a PV number then the variational principle for Hausdorff dimension does not hold 
for $([-1,1]^{2},f_{\beta})$ i.e. $\{\dim_{H}\mu|\mu ~f_{\beta}\mbox{-ergodic}\}<2$.
\end{theorem}
\begin{remark}
Theorem 2.1 is an extension of the result of Alexander and York \cite{[AY]}
that states that the Sinai-Ruelle-Bowen measures for 
$([-1,1]^{2},f_{\beta})$ does not have full R\'enyi dimension.
\end{remark}
\begin{remark}
It follows from \cite{[AY]} together with Solomyak's theorem 
about Bernoulli convolutions \cite{[SO1]} that for almost all
$\beta\in (0.5,1)$ the Sinai-Ruelle-Bowen measures for 
$([-1,1]^{2},f_{\beta})$ has full dimension. Thus our theorem shows
that in situations where there generically is an ergodic measure of full dimension the variational principle for Hausdorff dimension may not hold in general
because of special number theoretical properties of some parameter values.
As far as we know our theorem provides the first example of this type.  
\end{remark}
Now we come to our second class of examples. For $\beta\in (0.5,1)$ and $\tau\in(0,0.5)$ we define two affine contractions
on $[-1,1]^{2}$ by  
\[ T_{1}^{\beta,\tau}(x,z)=(\beta x+(1-\beta),\tau z+(1-\tau)) \]
\[ T^{\beta,\tau}_{-1}(x,z)=(\beta x-(1-\beta),\tau z-(1-\tau)). \]
From $\cite{[HU]}$ we know that there is a unique compact self-affine subset $\Lambda_{\beta,\tau}$ of
$[-1,1]^{2}$ satisfying \[\Lambda_{\beta,\tau}=T_{1}^{\beta,\tau}(\Lambda_{\beta,\tau})\cup
T^{\beta,\tau}_{-1}(\Lambda_{\beta,\tau}).\] 
Let $T_{\beta,\tau}$ be the smooth
expanding transformation 
on $T_{1}^{\beta,\tau}([-1,1]^{2})\cup T_{-1}^{\beta,\tau}([-1,1]^{2})$
defined by
\[T_{\beta,\tau}(x)=(T_{i}^{\beta,\tau})^{-1}(x)~~ \mbox{ if }~~ x\in T_{i}^{\beta,\tau}([-1,1]^{2})~~ \mbox{ for } ~~~i=1,-1.\]
Obviously the set $\Lambda_{\beta,\tau}$ is an invariant repeller for the transformation $T_{\beta,\tau}$. 
We call the system $(\Lambda_{\beta,\tau},T_{\beta,\tau})$ a
{\bf self-affine repeller}.\newline\newline 
Let us state our main result about the systems $(\Lambda_{\beta,\tau},T_{\beta,\tau})$.
\begin{theorem}
Let $\beta\in (0.5,1)$ be the reciprocal of a PV number. For all $\tau\in (0,0.5)$ we have 
$\dim_{H}\Lambda_{\beta,\tau}<\dim_{B}\Lambda_{\beta,\tau}$. Moreover if $\tau$ is sufficient small there
can not be a Bernoulli measure of full dimension for the system $(\Lambda_{\beta,\tau},T_{\beta,\tau})$.
\end{theorem}
\begin{remark}
We know from \cite{[NE]} that for almost all $\beta\in(0,5,1)$ and all $\tau\in(0,0.5)$ the identity
\[ \dim_{H}\Lambda_{\beta,\tau}=\dim_{B}\Lambda_{\beta,\tau}=\frac{\log2\beta}{\log\tau}+1 \]
holds and that there is a Bernoulli measure of full dimension
for $(\Lambda_{\beta,\tau},T_{\beta,\tau})$. Thus Theorem 2.2 shows that
dimensional theoretical properties
of dynamical systems can considerably change because of number theoretical
peculiarities.
\end{remark}
\begin{remark}
That the identity of Hausdorff and box-counting dimension may drop because of number
theoretical peculiarities was shown before by Przytycki and Urbanski \cite{[PU]}
in the context of Weierstrass like functions.  Pollicott and Wise \cite{[PW]} claimed (without a proof) that the first statement
of our theorem follows for small $\tau$ from the work of Przytycki and Urbanski. We were not able
to see that this is true and thus wrote down an independent proof which gives explicit upper bounds on $\dim_{H}\Lambda_{\beta,\tau}$ (see section seven). 
\end{remark} 
\begin{remark}
We do not know if there exists an ergodic measure of full Hausdorff dimension
for the systems $(\Lambda_{\beta,\tau},T_{\beta,\tau})$ and we can not calculate $\dim_{H}\Lambda_{\beta,\tau}$
in the case that $\beta\in(0.5,1)$ is the reciprocal of a PV number.
The second statement of our theorem only shows that it is not possible to calculate $\dim_{H}\Lambda_{\beta,\tau}$ by means of
Bernoulli measures in this situation. 
\end{remark}
\section{Coding maps and representation of ergodic measures}
We first introduce here the symbolic spaces which we use for
our coding. Let $\Sigma=\{-1,1\}^{\Z}$ and $\Sigma^{+}=\{-1,1\}^{\N_{0}}$. 
By $pr_{+}$ we denote the projection from $\Sigma$ onto $\Sigma^{+}$.
With a natural 
product metric 
$\Sigma$ (resp. $\Sigma^{+}$) becomes a perfect, totally disconnected and compact metric space. 
For $u,v\in \Z$ (resp. $u,v\in \N$) and $t_{0},t_{1}\dots,t_{u}\in\{-1,1\}$ we define a cylinder set in $\Sigma$ (resp. $\Sigma^{+}$) by
\[[t_{0},t_{1}\dots,t_{u}]_{v}:=\{(s_{k})|s_{v+k}=t_{k}
\mbox{ for } k=0,\dots, u\}.\]
The cylinder sets form a basis for the metric topology on $\Sigma$ (resp. $\Sigma^{+}$).
The forward shift map $\sigma$ on $\Sigma$ (resp. $\Sigma^{+}$) is given by
$\sigma((s_{k}))=(s_{k+1})$. The backward shift
$\sigma^{-1}$ is defined on $\Sigma$ and given by 
$\sigma((s_{k}))=(s_{k-1})$. By $b^{p}$ for $p\in(0,1)$ we denote the Bernoulli measure on
$\Sigma$ (resp. $\Sigma^{+}$), which is the product of the discrete
measure giving $1$ the probability $p$ and $-1$ the
probability $(1-p)$. We write $b$ for the equal-weighted 
Bernoulli measure $b^{0.5}$.
The Bernoulli measures are ergodic with respect to forward and
backward shifts (see \cite{[DGS]}).\newline\newline
We are now prepared to define the Shift coding for the Fat Baker's transformation $([-1,1]^{2},f_{\beta})$. Define a continuous map $\hat\pi_{\beta}$ from $\Sigma$ onto $[-1,1]^{2}$ by
\[ \hat\pi_{\beta}(\underline{i})=((1-\beta)\sum_{k=0}^{\infty}i_{k}\beta^{k},\sum_{k=1}^{\infty}i_{-k}(1/2)^{k}).\]
A simple check shows that
\[ f_{\beta}\circ\hat\pi_{\beta}(\underline{i})=\hat\pi_{\beta}\circ\sigma^{-1}(\underline{i}) \quad \forall\underline{i}\in\bar\Sigma=(\Sigma\backslash\{(s_{k})| 
\exists k_{0}  
\forall k \le k_{0} :s_{k}=1 \rbrace)\cup\{(1)\}.\]
Note that if $\mu$ is a $\sigma$-invariant Borel probability measure on $\Sigma$ we have $\mu(\bar\Sigma)=1$.
From this fact 
by applying standard techniques in ergodic theory it is possible to show that the
map
\[ \mu\longmapsto \mu_{\beta}:=\mu\circ \hat\pi_{\beta}^{-1}\]
from the space of $\sigma$-ergodic Borel probability measures on $\Sigma$ 
is continuous with respect to the weak$^{*}$ topology
and is onto the space of $f_{\beta}$-ergodic Borel probability measures on $[-1,1]^{2}$.
Moreover the system $([-1,1],f_{\beta},\mu_{\beta})$ is a measure theoretical factor of $(\Sigma,\sigma^{-1},\mu)$
\newline\newline
Now we introduce a shift coding for the self-affine repeller
$(\Lambda_{\beta,\tau},T_{\beta,\tau})$. Consider the homeomorphism $\pi_{\beta,\tau}\colon \Sigma^{+}\to \Lambda_{\beta,\tau}$ given by
\[ \pi_{\beta,\tau}(\underline{i})=((1-\beta)\sum_{k=0}^{\infty}i_{k}\beta^{k},(1-\tau)\sum_{k=0}^{\infty}i_{k}\tau_{k}).\]
It is easy to see that $\pi_{\beta,\tau}\circ\sigma=T_{\beta,\tau}\circ\sigma$.
Thus the systems $(\Lambda_{\beta,\tau},T_{\beta,\tau})$ is homoeomorph conjugated to $(\Sigma,\sigma)$ and the map
\[ \mu\longmapsto \mu_{\beta,\tau}:=\mu\circ\pi_{\beta,\tau}^{-1} \]
is a homeomorphism with respect to the weak$^{*}$ from the space of $\sigma$-ergodic Borel probability measures on $\Sigma^{+}$ onto the space of $T_{\beta,\tau}$-ergodic Borel probability measures on $\Lambda_{\beta,tau}$.
\section{Erd\"os measures and Garsia entropy}
For $\beta\in(0.5,1)$ define a continuous map from $\Sigma^{+}$ onto 
$[-1,1]$ by
\[ \pi_{\beta}(\underline{i})=(1-\beta)\sum_{k=0}^{\infty}i_{k}\beta^{k}. \]
Given a Borel probability measure $\nu$ on $\Sigma^{+}$ we define a Borel probability measure
on $[-1,1]$ by $\nu_{\beta}=\nu\circ\pi_{\beta}^{-1}$. If we choose the Bernoulli measure
$b^{p}$ on $\Sigma^{+}$ for a $p\in(0,1)$ then $b^{p}_{\beta}$ is a self-similar measure which is usually a called Bernoulli convolution.
There are a lot of results in the literature about Bernoulli convolutions and we can not cite all these works here.  Instead
we like to refer to the nice overview article ''Sixty years of Bernoulli convolutions''
by Peres, Schlag and Solomyak \cite{[PSS]}. \newline\newline
In our work we are not only interested in Bernoulli convolutions
but in all measures $\nu_{\beta}$ where $\nu$ is a $\sigma$ invariant Borel probability measure $\Sigma^{+}$
and $\beta\in(0.5,1)$ is the reciprocal of a PV number (see appendix B). We call a measure
of this type an {\bf Erd\"os} measure. \newline\newline 
Now we will introduce a special kind of entropy
related to Erd\"os measure. What we will do here is generalisation of the approach of Garsia (\cite{[GA1]},\cite{[GA2]}) for Bernoulli convolutions to all Erd\"os measures.
Let $\sim_{n,\beta}$ be the equivalence relation on $\Sigma^{+}$
given by
\[ \underline{i}~\sim_{n,\beta}~\underline{j}~\Leftrightarrow~\sum_{k=0}^{n-1}i_{k}\beta^{k}~=~\sum_{k=0}^{n-1}j_{k}\beta^{k} \]
and define a partition $\Pi_{n,\beta}$ of $\Sigma^{+}$ by $\Pi_{n,\beta}=\Sigma^{+}/\sim_{n,\beta}$. 
Recall that entropy of
a partition $\Pi$ with respect to a Borel probability measure $\nu$ on $\Sigma^{+}$ is
\[ H_{\nu}(\Pi)=-\sum_{P\in\Pi}\nu(P)\log\nu(P).\]
We denote the join of two partitions $\Pi_{1}$ and
$\Pi_{2}$ by $\Pi_{1}\vee\Pi_{2}$. This is the partition consisting of all sections $A\cap B$
for $A\in\Pi_{1}$ and $B\in\Pi_{2}$. It is easy to see that the 
 $\Pi_{n,\beta}\vee
\sigma^{-n}(\Pi_{m,\beta})$
is finer than the partition $\Pi_{n+m,\beta}$ and hence the sequence $H_{\nu}(\Pi_{n,\beta})$ 
is sub-additive for a shift invariant measure $\nu$ on $\Sigma^{+}$.
We can thus define the {\bf Garsia entropy} $G_{\beta}(\nu)$ for
a shift invariant Borel probability measure $\nu$ on $\Sigma^{+}$ by
\[ G_{\beta}(\nu):=\lim_{n\longrightarrow\infty}\frac{H_{\nu}(\Pi_{n,\beta})}{n}= \inf_{n}\frac{H_{\nu}(\Pi_{n,\beta})}{n}.\] 
The limit exists and is equal to the infimum since the sequence $H_{\nu}(\Pi_{n,\beta})$ is sub-additive. Another simple consequence of the sub-additivity of this sequence is that the map
\[ \nu\longmapsto G_{\beta}(\nu) \]
upper-semi-continuous with respect to the weak$^{*}$ topology on the space
of $\sigma$ invariant Borel probability measures on $\Sigma^{+}$.\newline\newline
We are now prepared to state our main theorem about Erd\"os measures and Garsia entropy.
\begin{theorem}
Let $\beta\in (0.5,1)$ be the reciprocal of a PV number. For all
$\sigma$-ergodic Borel probability measures $\nu$ on $\Sigma^{+}$ the following equivalence holds
\[ \nu_{\beta}\mbox{ is singular}\Leftrightarrow G_{\beta}(\nu)<-\log\beta     \Leftrightarrow \dim_{H}\nu_{\beta}<1.\]
Moreover the set of $\sigma$-ergodic measures Borel probability measures $\nu$ on $\Sigma^{+}$ such that $\nu_{\beta}$ is singular 
is open in the weak$^{*}$ topology and contains the Bernoulli measures $b^{p}$ for $p\in(0,1)$.
\end{theorem}
\begin{remark}
It has been shown by Erd\"os \cite{[ER1]} that the equal-weighted Bernoulli convolution $b_{\beta}$ is
singular if $\beta\in (0.5,1)$ is the reciprocal of a PV number. Using this result Garsia \cite{[GA2]} proved
$G_{\beta}(b)<-\log\beta$ and from this Alexander and Yorke \cite{[AY]} deduced
that the R\'enji dimension of $b_{\beta}$ is less than one. In the proof
of Theorem 4.1 we will adopt ideas of all of these authors. 
In our generalisation from the equal-weighted
Bernoulli measure to all $\sigma$-invariant measures
we had do deal with some difficulties which are mainly of technical
nature (see section four). 
\end{remark}
\begin{remark}
The PV case is exceptional. It was shown by Solomyak \cite{[SO1]} that for almost all
$\beta\in(0.5,1)$ the Bernoulli convolution $b_{\beta}$ is absolutely continuous with density in $L^{2}$.
\end{remark}
\section{Proof of Theorem 4.1}
The proof of Theorem 4.1 follows from three propositions and is given at the end of this section
\begin{proposition}
If $\beta\in(0,5,1)$ is
the reciprocal
PV number then the measures $b^{p}_{\beta}$ are singular
for all $p\in(0,1)$.
\end{proposition} 
\proof
The measure $b^{p}_{\beta}$ is given by the infinite convolution
of the discrete measures $b^{p,n}_{\beta}$, which give $(1-\beta)\beta^{n}$ the probability $p$ and $-(1-\beta)\beta^{n}$
the probability $(1-p)$. 
From \cite{[JW]} we know that
the Fourier transformation of a convolution is the product of
the Fourier transformation of the convolved measures. 
Consequently the Fourier transformation $\phi$ 
of $b^{p}_{\beta}$ is given by:
\[ \phi(b^{p}_{\beta},\omega)=\prod_{n=0}^{\infty}(\cos((1-\beta)\beta^{n}\omega)+(2p-1)\sin((1-\beta)\beta^{n}\omega)). \]
We see that
\[|\phi(b^{p}_{\beta},\omega)|=\prod_{n=0}^{\infty}|(\cos((1-\beta)\beta^{n}\omega)+(2p-1)\sin((1-\beta)\beta^{n}\omega))|\]
\[\ge
\prod_{n=0}^{\infty}|\cos((1-\beta)\beta^{n}\omega)|.\]
Now let $\omega_{k}=2\pi\beta^{-k}/(1-\beta)$. We have
\[ |\phi(b^{p}_{\beta},\omega_{k})|\ge\prod_{n=0}^{\infty}|\cos(2\pi\beta^{n-k})|=\prod_{n=0}^{k}|\cos(2\pi\beta^{n-k})|\prod_{n=k+1}^{\infty}|\cos(2\pi\beta^{n-k})|\]\[
=
C\prod_{n=0}^{k}|\cos(2\pi\beta^{-n})|\]
where $C$ is a constant independent of $k$ and not zero. Now let $\beta$ be the reciprocal of a PV number.
From proposition B1 of appendix B we know 
that there is a constant $0<\theta<1$ such that
$||\beta^{-n}||_{ \Z}\le\theta^{n}~ \forall n\ge0$ where 
$||.||_{ \Z}$ denotes the distance to the nearest integer. 
This implies $|\phi(b^{p}_{\beta},\omega_{k})|\ge \hat C>0$ for all $k>0$.
Thus we have that $|\phi(b^{p}_{\beta},\omega)|$ does not tend to zero
with $\omega\longrightarrow\infty$. 
Hence by Riemann-Lebesgue lemma $b^{p}_{\beta}$ can not be 
absolutely continuous if $\beta$ is the reciprocal of a PV number. But it follows from the theory of infinity
convolutions developed by Jessen and Winter \cite{[JW]} that $b^{p}_{\beta}$ is 
of pure type that means either absolutely continuous.  This completes the proof.  \qed
\begin{remark}
This proof is nothing but an obvious extension of Erd\"os \cite{[ER1]} original argument.
\end{remark}
\begin{proposition}
Let $\beta\in(0.5,1)$ be the reciprocal of a PV number and $\nu$ be a shift invariant Borel probability measure on $\Sigma^{+}$. 
If
$\nu_{\beta}$ is singular then $G_{\beta}(\nu)<-\log\beta$ holds.
\end{proposition}
\proof
Fix $\beta$. Define $\pi_{n}$ from $\Sigma^{+}$ to $[-1,1]$ by
$\pi_{n}((s_{k}))=\sum_{k=0}^{n-1}s_{k}(1-\beta)\beta^{k}$ and let $\nu_{n}=\nu\circ\pi_{n}^{-1}$.
Let $\sharp(n)$ be the number of distinct points of the form $\sum_{k=0}^{n-1}\pm(1-\beta)\beta^{k}$
and $\omega(n)$ be the minimal distance between two of those points. 
Furthermore denote the points by $x_{i}^{n}~~i=1\dots\sharp(n)$ and
let $m_{i}^{n}$ be the $\nu$ measure
of the corresponding elements in $\Pi_{n,\beta}$, which means $m_{i}^{n}=\nu_{n}(x_{i}^{n})$. \newline\newline
We first state a property of PV numbers we will have to use here, 
see proposition B2 of appendix B:
\[ \beta^{-1} \mbox{ is PV number }\Rightarrow \exists ~~\bar c:\omega(n)\ge \bar c\beta^{n}.     \]
Since $(\sharp(n)-1)\omega(n)\le 2$ we get $\sharp(n)\le 4\omega(n)^{-1}\le c\beta^{-n}$ with $c:=4\bar c^{-1}$.\newline\newline
Now we assume that
$\nu_{\beta}$ is singular.  It follows that there exists a constant $C$ such that: \newline
$\forall \epsilon >0~\exists $ disjoint intervals $(a_{1},b_{1}),\dots,(a_{u},b_{u})$ with 
\[ \sum_{l=1}^{u}(b_{l}-a_{l}) < \epsilon~~ \mbox{ and }~~ \nu_{\beta}(O)>C
~~\mbox{where}~~O:=\bigcup_{l=0}^{u}(a_{l},b_{l}).\]
With out loss of generality we may assume 
$\nu_{\beta}(a_{l})=\nu_{\beta}(b_{l})=0$ for $l=1\dots u$.
It is obvious that the discreet distribution $\nu_{n}$ converges  weakly to $\nu_{\beta}$. 
Thus we have: $\exists n_{1}(\epsilon)~~\forall n>n_{1}(\epsilon)~:~\nu_{n}(O)>C $.
We now expand the intervals a little bit, so that their
length is a multiple of $\omega(n)$.
\[ k_{l,n}:=\max\{k~|~ k \omega(n) \le a_{l} \} \quad a_{l,n}:=k_{l,n}\omega(n) \] 
\[ \bar k_{l,n}:=\min\{k~|~ b_{l} \le k \omega(n) \} \quad b_{l,n}:=\bar k_{l,n}\omega(n) \]
Since $\omega(n)\longrightarrow 0$ we have: \newline
$\exists n_{2}(\epsilon)>n_{1}(\epsilon)~~\forall n>n_{2}(\epsilon):
~~(a_{l,n},b_{l,n})$ disjunct for $l=1 \dots u$ and
\[ \sum_{l=1}^{u}(b_{l,n}-a_{l,n}) < \epsilon~~ \mbox{ and }~~ \nu_{n}(\bar O)>C
~~\mbox{ where }~~\bar O=\bigcup_{l=0}^{u}(a_{l,n},b_{l,n}).\]
Let $\hat \sharp(n)$ be the number of distinct points $x_{i}^{n}$ in $\bar O$. Since in one
interval $(a_{l,n},b_{l,n})$ there are at most $\bar k_{l,n}-k_{l,n}$ points
$x_{i}^{n}$ we have $\omega(n)\hat\sharp(n)\le\epsilon$ and
hence $\hat\sharp(n)\le\epsilon c \beta^{-n}$.  \newline
For all $n>n_{2}(\epsilon)$ we can now estimate:
\[ H_{\nu}(\Pi_{n,\beta})=-\sum_{i=1}^{\sharp(n)}m_{i}^{n}\log m_{i}^{n}=- 
\sum_{x_{i}^{n}\in\bar O}m^{n}_{i}\log m^{n}_{i}-\sum_{x_{i}^{n}\not\in\bar O}m^{n}_{i}\log m^{n}_{i}\] \[
\le \nu_{n}(\bar O)\log\frac{\hat\sharp(n)}{\nu_{n}(\bar O)}+(1-\nu_{n}(\bar O))\log\frac{\sharp(n)-\hat\sharp(n)}{1-\nu_{n}(\bar O)}
\] \[\le \nu_{n}(\bar O)\log\hat\sharp(n)+(1-\nu_{n}(\bar O))\log\sharp(n)+\log 2\]
\[  \le \nu_{n}(\bar O)\log\epsilon c\beta^{-n}+(1-\nu_{n}(\bar O))\log c\beta^{-n}+\log 2 \] \[
\le n\log\beta^{-1}+C\log\epsilon+\log c +\log 2 .\]
If $\epsilon$ is small enough we have $H_{\nu}(\Pi_{n,\beta})/n<\log\beta^{-1}$ for all $n\ge n_{2}(\epsilon)$.
Using the sub-additivity of $H_{\nu}(\Pi_{n,\beta})$ we get our result.\qed
\begin{remark}
Garsia sketched a proof of this proposition for the equal weighted Bernoulli measure
in \cite{[GA1]}
. Our poof is a more detailed
and extended version of Garsia's argumentation.
\end{remark} 
\begin{proposition}
If $\nu$ is a shift ergodic Borel probability measure on $\Sigma^{+}$ and $\beta\in(0.5,1)$ we have 
\[\dim_{H}\nu_{\beta}\le G_{\beta}(\nu)/-\log\beta.\]
\end{proposition}
\proof Because we will operate with R\'enyi dimension $\dim_{R}$ (see appendix A)
we are interested in an upper bound on the quantity
\[ h_{\nu}(\epsilon)=\inf\{H_{\nu}(\Pi)|\Pi \mbox{ a partition with } \mbox{diam}\Pi\le\epsilon\}\] by the entropy
of the partitions $\Pi_{n,\beta}$ of $\Sigma^{+}$. We proof the following statement

\[ h_{\nu_{\beta}}(2\beta^{n})\le H_{\nu}(\Pi_{n,\beta}).\]
Fix $\beta\in(0.5,1)$, $\tau\in(0,0.5)$, a measure $\nu$ on $\Sigma^{+}$
and $n\in\N$. We use the convention that the first coordinate axis
is called $x$-axis and $pr_{X}$ denotes the projection on this axis.\newline
We define a partition of $\Lambda_{\beta,\tau}$ by $\wp_{n}=\pi_{\beta,\tau}(\Pi_{n,\beta})$.
By definition we have \[ H_{\nu}(\Pi_{n,\beta})=H_{\nu_{\beta,\tau}}(\wp_{n}). \]
We should
say something about the structure of $\wp_{n}$.
The image of a cylinder set $[i_{0},\dots,i_{n-1}]_{0}$ in $\Sigma^{+}$ under $\pi_{\beta,\tau}$
is the part of $\Lambda_{\beta,\tau}$ lying in the rectangle $T_{i_{n-1}}^{\beta,\tau}\circ\dots\circ T_{i_{0}}^{\beta,\tau}(Q)$
of $x$-length $2\beta^{n}$. It is not difficult to check
that two cylinder sets lie in the same element
of $\Pi_{n,\beta}$ 
if and only if the corresponding rectangles lie above each
other. So the projection of
an element in $\wp_{n}$ onto the $x$-axis has length $2\beta^{n}$.\newline
The projection onto the $x$-axis of two elements in $\wp_{n}$ may overlap.
Starting with $\wp_{n}$, we want to construct inductively a partition
$\bar\wp_{n}$ of $\Lambda_{\beta,\tau}$ with
non-overlapping projections,
in a way that does neither increase length of the
projections nor entropy. Let $N(\wp)$ be the number of pairs of
elements in a partition $\wp$ that do have overlapping projections onto the
$x$-axis. We now construct a finite sequence $\wp_{n}^{k}$ of partitions.
First let $\wp_{n}^{0}=\wp_{n}$. Now let $\wp_{n}^{k}$ be constructed
and $N(\wp_{n}^{k})>0$.\
Let $P_{1}$ and $P_{2}$ be two elements of
$\wp_{n}^{k}$ with overlapping projections. Without loss of generality we may assume $\nu_{\beta,\tau}(P_{1})\ge\nu_{\beta,\tau}(P_{2})$
and define:
\[\hat P_{1}=P_{1}\cup(P_{2}\cap(pr_{X}P_{1}\times [-1,1]))\qquad
\hat P_{2}=P_{2}\backslash (pr_{X}P_{1}\times [-1,1]) .    \]
We have $\hat P_{1}\dot\cup \hat P_{2}=P_{1}\dot\cup P_{2}$,
$P_{1}\subseteq\hat P_{1}$ and $\hat P_{2}\subseteq P_{2}$.
Thus we know:
$\nu_{\beta,\tau}(P_{1})+\nu_{\beta,\tau}(P_{2})=\nu_{\beta,\tau}(\hat P_{1})+\nu_{\beta,\tau}(\hat P_{2})$ and
$\nu_{\beta,\tau}(\hat P_{1})\ge\nu_{\beta,\tau}(P_{1})
\ge \nu_{\beta,\tau}(P_{2}) \ge \nu_{\beta,\tau}(\hat P_{2})$. 
Since the function $-x\log x$ is concave, this implies:
\[-(\nu_{\beta,\tau}(\hat P_{1})\log\nu_{\beta,\tau}(\hat P_{1})+
   \nu_{\beta,\tau}(\hat P_{2})\log\nu_{\beta,\tau}(\hat P_{2})) \le \]
\[-(\nu_{\beta,\tau}(P_{1})\log\nu_{\beta,\tau}(P_{1})+
   \nu_{\beta,\tau}(P_{2})\log\nu_{\beta,\tau}(P_{2})). \]
Hence if we substitute $\hat P_{1}$, $\hat P_{2}$ for $P_{1}$, $P_{2}$, we get
a partition $\wp^{k+1}_{n}$ of $\Lambda_{\beta,\tau}$
with non-increased entropy. 
From the definition of $\hat P_{1}$ and $\hat P_{2}$ we see
that $pr_{X}\hat P_{1}=pr_{X}P_{1}$, $pr_{X}\hat P_{2}\subseteq pr_{X}P_{2}$
and that the projections of $\hat P_{1}$ and $\hat P_{2}$ onto the
$x$-axis do not overlap. So the length of the projections are obviously not increased. Furthermore we observe
that there cannot be any new overlaps of the projections of $\hat P_{1}$ or $\hat P_{2}$ with
the projections of other elements in $\wp^{k}_{n}$, that do
not appear, when we consider $P_{1}$ or $P_{2}$.
Hence
$N(\wp_{n}^{k+1})<N(\wp_{n}^{k})$. \newline\newline
So after a finite number of steps we get
a partition $\bar\wp_{n}$ with
\[ H_{\nu_{\beta,\tau}}(\wp_{n})\ge H_{\nu_{\beta,\tau}}(\bar\wp_{n}), \]
non-overlapping projections onto the $x$-axis and $\mbox{diam}~pr_{X}\bar\wp_{n}\le2\beta^{n}$. 
$pr_{X}\bar\wp_{n}$ is a partition of the interval $[-1,1]$ and we
have
\[ H_{\nu_{\beta}}(pr_{X}\bar\wp_{n})=H_{\nu_{\beta,\tau}}(\bar\wp_{n}), \]
since the measure $\nu_{\beta}$ is the projection of $\nu_{\beta,\tau}$
onto the $x$-axis. The proof of our claim is complete:
\[ h_{\nu_{\beta}}(2\beta^{n})\le
H_{\nu_{\beta}}(pr_{X}\bar\wp_{n})=H_{\nu_{\beta,\tau}}(\bar\wp_{n})\le
H_{\nu_{\beta,\tau}}(\wp_{n})=
H_{\nu}(\Pi_{n,\beta}).  \]
We are now able to estimate the R\'enyi dimension
\[ \overline{\dim}_{R}\nu_{\beta}
=\overline{\lim}_{\epsilon \longrightarrow \infty}\frac{h_{\nu_{\beta}}(\epsilon)}{\log\epsilon^{-1}}
=\overline{\lim}_{n \longrightarrow \infty}
\frac{h_{\nu_{\beta}}(2\beta^{n})}{\log 0.5\beta^{-n}}=
\overline{\lim}_{n \longrightarrow \infty}
\frac{h_{\nu_{\beta}}(2\beta^{n})}{n\log \beta^{-1}}\]
\[
\le \lim_{n \longrightarrow \infty}\frac{H_{\nu}(\Pi_{n,\beta})}{n\log \beta^{-1}}=\frac{G_{\beta}(\nu)}{\log\beta^{-1}} .\]     \newline
Using part (3) of proposition A1 from appendix A we get \newline
\[ \forall \delta>0\exists X:\nu_{\beta}(X)>0
\mbox{ and } \underline{d}(x,\nu_{\beta})\le G_{\beta}(\nu)/\log\beta^{-1}+\delta\quad \forall x\in X.\] 
But the measure $\nu_{\beta}$ is exact dimensional, because
it is the transversal measure in the context of
the ergodic dynamical system $(\Lambda_{\beta,\tau},T_{\beta,\tau},\nu_{\beta,\tau})$. This fact was observed by Ledrappier and Porzio, see \cite{[LP]}.
So our estimate must hold
$\nu_{\beta}$-almost everywhere
and by part (2) of proposition A2 we get
$\dim_{H}\nu_{\beta}\le G_{\beta}(\nu)/\log\beta^{-1}+\delta$ for all $\delta>0$.
This proves the proposition.\qed
\begin{remark}
Let us remark that Alexander and Yorke \cite{[AY]} proved the 
identity $\dim_{R}b_{\beta}=G_{\beta}(b)/\log\beta^{-1}$
for the equal-weighted infinitely convolved Bernoulli measure $b_{\beta}$.
In their proof
they used the self-similarity 
of this measure. In our general situation we could 
not appeal to self-similarity and thus had to develop a
different technique.
\end{remark}
{\bf Proof of Theorem 4.1}
Under the assumptions of our theorem we have
\[ \nu_{\beta} \mbox{ is singular}\Rightarrow^{5.2}
G_{\beta}(\nu)<\log\beta^{-1}
\Rightarrow^{5.3}\dim_{H}\nu_{\beta}<1\Rightarrow\nu_{\beta} \mbox{ is singular} . \]
These implications prove the first statement of Theorem 4.1.
Now choose an singular Erd\"os measure $\xi_{\beta}$.
We have $G_{\beta}(\xi)<\log\beta^{-1}$. By upper-semi-continuity
of $G$ we get $G_{\beta}(\nu)<\log\beta^{-1}$ and hence
$\dim\nu_{\beta}<1$ for all $\nu$ in a hole $\mbox{weak}^{*}$ 
neighbourhood of $\xi$. Thus
the set $\{\nu|\nu_{\beta} \mbox{ is singular}\}$ is open
in the $\mbox{weak}^{*}$ topology. The set contains 
all Bernoulli measure by Proposition 4.1.\qed
\section{Proof of Theorem 2.1}
The proof of Theorem 2.1 follows from 
and Theorem 4.1
and two propositions
providing upper estimates on the Hausdorff dimension of all ergodic measures $\mu_{\beta}$ 
for the Fat Baker's transformation $f_{\beta}$. It can be found at the end of this section.
\begin{proposition}
If $\mu$ is a shift ergodic Borel probability measure on $\Sigma$ and $\beta\in(0,5)$ we have 
\[\dim_{H}\mu_{\beta}\le 1+\le G_{\beta}(pr_{+}(\mu))/-\log\beta\]
where $pr_{+}$ denotes the projection from $\Sigma$ onto $\Sigma^{+}$.
\end{proposition}
\proof By Proposition A2 and the definition of the Hausdorff dimension of a measure we have $\dim_{H}\mu_{\beta}\le 1+\dim_{H}pr_{X}\mu_{\beta}$ where $pr_{X}$ denotes the projection onto the first coordinate axis. Just by definition of the involved measures we have
$pr_{X}\mu_{\beta}=(pr_{+}\mu)_{\beta}$ and hence $\dim_{H}\mu_{\beta}\le 1+\dim_{H}(pr_{+}\mu)_{\beta}$. The proposition follows
now immediately from Proposition 5.3.\qed
\begin{proposition}
If $\mu$ is a shift ergodic Borel probability measure on $\Sigma$ and $\beta\in(0,5)$ we have 
\[\dim_{H}\mu_{\beta}\le 1+\le h_{\mu}(\sigma)/\log 2\]
where $h_{\mu}(\sigma)$ is the usual measure-theoretic entropy of the shift $(\Sigma,\sigma,\mu)$.
\end{proposition}
\proof The proof of this proposition is a little bit difficult. We want to use the general theory
relating the dimension of ergodic measure to entropy and Lyapunov exponents (see \cite{[LY]} and \cite{[BPS]}). Usually
this theory is stated in the context of diffeomorphisms but the Fat Baker's transformation is not invertible and has
a singularity. To deal with the first problem we define for $\beta\in(0.5,1)$ and $\tau \in (0,0.5)$
a lift $\hat f_{\beta,\tau}\colon [-1,1]^{3}\to [-1,1]^{3}$ of the Fat Baker's transformation $f_{\beta}$ by
\[  
{\hat f_{\beta} (x,y,z) } = \lbrace
\begin{array}{cc}
(\beta x+(1-\beta),2y-1,\tau z+(1-\tau))\quad \mbox{if}\quad y \geq 0  \\
(\beta x-(1-\beta),2y+1,\tau z-(1-\tau))\quad \mbox{if} \quad y<0    .
\end{array}
\]
This maps is invertible and its projection onto the $(x,y)$-plane is $f_{\beta}$. Moreover it is easy to
see that $f_{\beta,\tau}$ has an attractor $\hat\Lambda_{\beta,\tau}$ which is given by the product
of the self-affine set $\Lambda_{\beta,\tau}$ in the $(x,z)$-plane with the interval $[-1,1]$ on the $y$-axis. Let us
introduce a Shift coding $\hat\pi_{\beta,\tau}:\Sigma\longmapsto\hat\Lambda_{\beta,\tau}$ for the system
$(\hat\Lambda_{\beta,\tau},\hat f_{\beta,\tau})$ by
\[ \hat\pi_{\beta,\tau}(\underline{i})=((1-\beta)\sum_{k=0}^{\infty}i_{k}\beta^{k},\sum_{k=1}^{\infty}i_{-k}(1/2)^{k},(1-\tau)\sum_{k=0}^{\infty}i_{k}\tau_{k}).\]
Given a $\sigma$-ergodic measure on $\Sigma$ we define a $\hat f_{\beta,\tau}$-ergodic measure $\hat\mu_{\beta,\tau}$ on $\hat\Lambda_{\beta,\tau}$ by $\hat\mu_{\beta,\tau}=\mu\circ\hat\pi_{\beta,\tau}^{-1}$. Section four of \cite{[NE]} contains a proof
of the fact that we are allowed to apply the general results found in \cite{[LY]} and \cite{[BPS]} to the system $(\hat\Lambda_{\beta,\tau},\hat f_{\beta,\tau},\hat \mu_{\beta,\tau})$ also this system has a singularity. We do not want to reproduce the argument here. We only like to mention that main idea is that the set of points that approaches the singularity of $(\hat\Lambda_{\beta,\tau},\hat f_{\beta,\tau},\hat \mu_{\beta,\tau})$ with
exponential speed has zero measure and thus Lyapunov charts exist almost everywhere for $(\hat\Lambda_{\beta,\tau},\hat f_{\beta,\tau},\hat \mu_{\beta,\tau})$. From Theorem C and Theorem F of \cite{[LY]} we have by this fact
\[ \dim_{H}\hat\mu_{\beta,\tau}\le \frac{h_{\hat \mu_{\beta,\tau}}(\hat f_{\beta,\tau})}{\log 2}+\dim\hat\mu_{\beta,\tau}^{s}.\]
where $\dim\hat\mu^{s}_{\beta,\tau}$ is the local dimension of the conditional measures of $\hat\mu_{\beta,\tau}$ on the
partition $\{[-1,1]\times \{y\} \times [-1,1]|y\in[-1,1]\}$ in the stable direction of $\hat f_{\beta,\tau}$ and $h_{\hat\mu_{\beta,\tau}}(\hat f_{\beta,\tau})$ is the measure theoretical entropy of the system $(\hat\Lambda_{\beta,\tau},\hat f_{\beta,\tau},\hat \mu_{\beta,\tau})$. Since the conditional measures are just by definition concentrated on
the set $\{(x,y,z)|(x,z)\in \Lambda_{\beta,\tau}~~y\in[-1,1]\}$ we have $\dim\mu_{\beta,\tau}^{s}\le \dim_{B}\Lambda_{\beta,\tau}$
and from \cite{[PW]} we know $\dim_{B}\Lambda_{\beta,\tau}=\log(2\beta/\tau)/\log(1/\tau)$.
Furthermore it is easy to see that the systems  $(\hat\Lambda_{\beta,\tau},\hat f_{\beta,\tau},\hat \mu_{\beta,\tau})$
and $(\Sigma,\sigma,\mu)$ are measure theoretical conjugated and thus $h_{\hat\mu_{\beta,\tau}}(\hat f_{\beta,\tau})=h_{\mu}(\sigma)$. 
Hence we have
\[ \dim_{H}\hat\mu_{\beta,\tau}\le \frac{h_{\mu}(\sigma)}{\log 2}+\frac{\log(2\beta/\tau)}{\log(1/\tau)}.\]
Now note that $\mu_{\beta,\tau}$ projects to $\mu_{\beta}$ and hence $\dim_{H}\mu_{\beta}\le\dim_{H}\hat\mu_{\beta,\tau}$
for all $\tau\in(0,0.5)$. Thus we get
\[ \dim_{H}\hat\mu_{\beta,\tau}\le \frac{h_{\mu}(\sigma)}{\log 2}+\frac{\log(2\beta/\tau)}{\log(1/\tau)}\forall \tau\in (0,0.5).\]
With $\tau\longrightarrow 0$ our proof is complete.\qed\newline\newline
{\bf Proof of Theorem 2.1}
From Theorem 4.1 and the upper-semi-continuity of $G_{\beta}$ we get $G_{\beta}(pr^{+}\mu)/\log\beta^{-1}\le c_{1}<1$ for all $\mu$ in hole $\mbox{weak}^{*}$ neighbourhood $U$ of $b$ in space of $\sigma$-ergodic Borel probability measures on $\Sigma$. 
Hence by Proposition 6.1 $\dim_{H}\bar\mu_{\beta}\le c_{1}+1<2$ holds for all $\mu$ in $U$.
On the other hand we have by well known properties of the measure theoretical entropy, $h_{\mu}(\sigma)/\log2\le c_{2}<1$ on the complement of $U$ (see \cite{[DGS]}).
From Proposition 6.1 we thus get $\dim_{H}\mu_{\beta}\le c_{2}+1<2$ for all $\mu$ in the complement of $U$. Putting these facts together
we obtain
\[ \dim_{H}\mu_{\beta}\le\max\{c_{1},c_{2}\}+1<2=\dim_{H}[-1,1]^{2}.\]
But we know that all ergodic measures for the system $([-1,1]^{2},f_{\beta})$ are of the form $\mu_{\beta}$ for
some $\sigma$-ergodic Borel probability measures $\mu$ on $\Sigma$. 
and the proof is complete.\qed
\section{Proof of Theorem 2.2}
The proof of Theorem 2.2 has a lot of ingredencies, a formula for $\dim_{B}\Lambda_{\beta,\tau}$ found in \cite{[PW]},
a formula for $\dim_{H}b^{p}_{\beta,\tau}$ found in \cite{[NE]}, Theorem 4.1 and the following two proposition giving upper bounds
on $\dim_{H}\Lambda_{\beta,\tau}$.
\begin{proposition}
If $\beta\in(0.5,1)$ is the reciprocal of an PV number and  $\tau \in (0,0.5)$ 
we have
\[ \dim_{H}\Lambda_{\beta,\tau}\le \frac{\log(\sum_{P\in\Pi_{n,\beta}}(\sharp P)^{
\frac{\log\beta}{\log\tau}})}{n\log\beta^{-1}}\qquad\forall n\ge 1\]
where $\Pi_{n,\beta}$ is the partition of $\Sigma^{+}$
defined in section four and $\sharp P$ denotes the number of
cylinder sets of length $n$ contained in an element of this partition.
\end{proposition}
\proof 
Fix a reciprocal of a PV number $\beta\in(0.5,1)$ and $\tau \in (0,0.5)$.
Let $n\ge 1$ and set \[u_{n}=\frac{\log(\sum_{P\in\Pi_{n,\beta}}(\sharp P)^{
\frac{\log\beta}{\log\tau}})}{n\log\beta^{-1}}.\]
Consider the set of cylinders in $\Sigma^{+}$
given by $C_{n}=\{[\tilde s_{1}\tilde s_{2}\dots\tilde s_{m}        ]_{0}~|~\tilde s_{i} \in\{-1,1\}^{n}~ i=1\dots m\}$.
Define a set function $\eta$ on $C_{n}$ by
\[ \eta([\tilde s]_{0})=\frac{\sharp P(\tilde s)^{\log\beta/\log\tau}}{\sharp P(\tilde s)}\beta^{nu_{n}}~\mbox{ and }~\]
\[\eta([\tilde s_{1}\tilde s_{2}\dots\tilde s_{m}]_{0})=\eta([\tilde s_{1}]_{0})\cdot\eta([\tilde s_{2}]_{0})\cdot\dots\cdot\eta([\tilde s_{m}]_{0})\]
where $\tilde s, \tilde s_{1},\dots \tilde s_{m}$ are elements of
$\{-1,1\}^{n}$ and $P(\tilde s)$ denotes the element of
the partition $\Pi_{n,\beta}$ containing the cylinder
$[\tilde s]_{0}$.\newline\newline
Note the facts that $C_{n}$ is a basis of the metric topology of $\Sigma^{+}$ and that $\sum_{\tilde s \in \{-1,1\}^{n}}\eta([\tilde s]_{0})=1$ by the definition of $u_{n}$. Thus we can extend $\eta$ to a Borel probability measure on $\Sigma^{+}$ 
and $\eta_{\beta,\tau}:=\eta\circ \pi_{\beta,\tau}^{-1}$ defines a Borel probability measure on $\Lambda_{\beta,\tau}$.\newline\newline
Given $m\ge1$ we set $q(m)=\lceil m(\log\beta/\log\tau)\rceil$. Given a $\tilde s_{i}\in \{-1,1\}^{n}$ for $i=1\dots m$ we define a subset of $\Lambda_{\beta,\tau}$ by
\[ R_{\tilde s_{1}\dots \tilde s_{n}}=
\{(\sum_{i=0}^{\infty} s_{i}(1-\beta)\beta^{i},\sum_{i=0}^{\infty}t_{i}(1-\tau)\tau^{i})~|~s_{i},t_{i}\in\{-1,1\} \] 
\[ (s_{(i-1)n},\dots,s_{in-1})=\tilde s_{i} \quad i=1\dots m\quad\mbox{and}\quad\] 
\[ (t_{(i-1)n},\dots,t_{in-1})=\tilde s_{i} \quad i=1\dots q(m) \}.\] 
We see that $R_{\tilde s_{1}\dots \tilde s_{m}}$ is "almost" a square in $\Lambda_{\beta,\tau}$ of side length $\beta^{mn}$. More precise we have:
\[ c_{1}\beta^{mn}\le \mbox{diam}  R_{\tilde s_{1}\dots \tilde s_{m}}\le c_{2} \beta^{mn} \qquad (1)\]
where the constants $c_{1},c_{2}$ are independent of the choice of
$\tilde s_{i}$.\newline\newline
Now let as examine the $\eta_{\beta,\tau}$ measure of the sets
$R_{\tilde s_{1}\dots \tilde s_{m}}$. \newline
Assume that
$\tilde t_{i}\sim_{n,\beta} \tilde s_{i}$ for $i=q(m)+1\dots m$ where $\sim_{n,\beta}$ is the equivalence relation introduced in section four.
The rectangles $\pi_{\beta,\tau}([\tilde s_{1}\dots \tilde s_{q(m)}\tilde t_{q(m)+1}\dots \tilde t_{m}]_{0})$ are all disjoint and lie above each other in the set $R_{\tilde s_{1}\dots \tilde s_{m}}$.
Hence we have
\[ \eta_{\beta,\tau}(R_{\tilde s_{1}\dots \tilde s_{m}})\ge \eta(\bigcup_{\tilde t_{i}\sim_{n,\beta} \tilde s_{i}~~ i=q(m)+1\dots m} \pi_{\beta,\tau}([\tilde s_{1}\dots \tilde s_{q(m)}\tilde t_{q(m)+1}\dots \tilde t_{m}]_{0}        )=\]
\[=\sum_{\tilde t_{i}\sim_{n,\beta} \tilde s_{i}~~ i=q(m)+1\dots m}\eta([\tilde s_{1}\dots \tilde s_{q(m)}\tilde t_{q(m)+1}\dots \tilde t_{m}]_{0}).\]
Using the fact $\tilde s\sim_{n,\beta} \tilde t\Rightarrow \sharp P(\tilde s)=\sharp P(\tilde t)\Rightarrow \eta([\tilde s]_{0})=\eta([\tilde t]_{0})$ this last expression equals
\[ \prod_{i=1}^{m}\eta([\tilde s_{i}]_{0})\sum_{\tilde t_{i}\sim_{n,\beta} \tilde s_{i}~~ i=q(m)+1\dots m}1\]
\[ =\prod_{i=1}^{m}
\frac{\sharp P(\tilde s_{i})^{\log\beta/\log\tau}}{\sharp P(\tilde s_{i})}\beta^{mnu_{n}}
\sum_{\tilde t_{i}\sim_{n,\beta} \tilde s_{i}~~ i=q(m)+1\dots m}1\]
\[
=\frac{\prod_{i=1}^{m}\sharp P(\tilde s_{i})^{\log\beta/\log\tau}}{\prod_{i=1}^{q(m)}\sharp P(\tilde s_{i})}  \beta^{ mnu_{n}}=(\phi_{\tilde s_{1}\dots \tilde s_{m}}\beta^{nu_{n}})^{m}\]
where
\[\phi_{\tilde s_{1}\dots \tilde s_{m}}
=(\frac{\prod_{i=1}^{m}\sharp P(\tilde s_{i})^{\log\beta/\log\tau}}{\prod_{i=1}^{q(m)}\sharp P(\tilde s_{i})})^{1/m} .\]
Now fix an $\epsilon>0$
We use the sets $R_{\tilde s_{1}\dots \tilde s_{m}}$ to construct a good cover of $\Lambda_{\beta,\tau}$ in the sense for Hausdorff dimension. To this end set 
\[R_{m}:=\{R_{\tilde s_{1}\dots \tilde s_{m}}|\phi_{\tilde s_{1}\dots \tilde s_{m}}\ge \beta^{n\epsilon} \}. \]
We have an upper bound on the cardinality of $R_{m}$.
If $R\in R_{m}$ then $\eta_{\beta,\tau}(R)\ge \beta^{mn(u_{n}+\epsilon)}$ and since $\eta_{\beta,\tau}$ is a probability measure we see:
\[ \mbox{card}(R_{m})\le \beta^{-mn(u_{n}+\epsilon)}\qquad(2).\]
~\newline
Now let $R(M)=\bigcup_{m\ge M} R_{m}$. We want to prove that
$R(M)$ is a cover of $\Lambda_{\beta,\tau}$ for all $M\ge 1$.\newline\newline
For $\underline s=(s_{k})\in \Sigma^{+}$ we define the function $\phi_{m}$ by
$\phi_{m}(\underline s)=\phi_{s_{0}\dots s_{mn-1}}$. In addition we need two auxiliary functions on $\Sigma^{+}$:
\[f_{m}(\underline{s})=\frac{\prod_{i=0}^{m} \sharp P((s_{(i-1)n},\dots,s_{in-1}))^{1/m}}{\prod_{i=0}^{q(m)}\sharp P((s_{(i-1)n},\dots,s_{in-1}))^{1/q(m)}},  \]
\[g_{m}(\underline s)=(\prod_{i=1}^{q(m)}\sharp P((s_{(i-1)n},\dots,s_{in-1})))^{1/q(m)(\log\beta\log\tau-q(m)/m)}.
\]
Since $1\le \sharp P(\tilde s)\le 2^{n}$ we have
$1\le g_{m}(\underline s)\le 2^{n(\log\beta/\log\tau-q(m)/m)}$.
Thus by the definition of $q(m)$ we have $g_{m}(\underline s)\longrightarrow 1$. Moreover we have $\overline{\lim}_{m\longrightarrow\infty}f_{m}(\underline{s})\ge 1$
because $\prod_{i=0}^{t} \sharp P((s_{i-1}n,\dots,s_{in-1}))^{1/t}\ge 1~~\forall t\ge 1$. \newline
A simple calculation shows $\phi_{m}(\underline{s})=(f_{m}(\underline{s}))^{\log\beta/\log\tau}g_{m}(\underline{s})$.
The properties of 
$f$ and $g$ thus imply:
\[ \overline{\lim}_{m\longrightarrow\infty}\phi_{m}(\underline{s})\ge1 \qquad \forall ~\underline{s}\in \Sigma^{+}. \]
This will help us to show that $R(M)$ is a cover of $\Lambda_{\beta,\tau}$. For all $\underline{s}=(s_{k})\in \Sigma^{+}$ there is an $m\ge M$ such that $\phi_{m}(\underline{s})\ge \beta^{n\epsilon}$ and thus $ 
\pi_{\beta,\tau}(\underline{s})\in R_{s_{0},\dots,s_{mn-1}}\in R(M)$.
Since $\pi_{\beta,\tau}$ is onto $\Lambda_{\beta,\tau}$ we see that 
$R(M)$ is indeed a cover of $\Lambda_{\beta,\tau}$.\newline\newline
We are now able to complete the proof. For every $\epsilon>0$
and every $M\in\N$ we have:
\[ \sum_{R\in R(M)}(\mbox{diam} R)^{u_{n}+2\epsilon}
=\sum_{m\ge M}\sum_{R\in R_{m}}(\mbox{diam} R)^{u_{n}+2\epsilon}\]
\[ \le^{(1)}\sum_{m\ge M}\sum_{R\in R_{m}}(c_{2}\beta^{mn})^{u_{n}+2\epsilon}=\sum_{m\ge M} \mbox{card}(R_{m})(c_{2}\beta^{mn})^{u_{n}+2\epsilon}\]
\[\le^{(2)} c_{2}^{u_{n}+2\epsilon}\sum_{m\ge M}\beta^{mn\epsilon}.\]
The last expression goes to zero with $M\longrightarrow 0$. By the definition for Hausdorff dimension we thus get $\dim_{H}\Lambda_{\beta,\tau}\le u_{n}+2\epsilon$ and since $\epsilon$ is arbitrary, we have $\dim_{H}\Lambda_{\beta,\tau}\le u_{n}$.\qed
\begin{remark}
Some ideas we used here are to due the prove of McMullen's theorem on self-affine carpets \cite{[MC]} by Pesin in \cite{[PE1]}.
\end{remark}
Now we use strategies developed in the proof of Proposition 5.2 to get:
\begin{proposition}
If $\beta\in(0.5,1)$ is the reciprocal of a PV number and  $\tau \in (0,0.5)$ 
we have
\[\exists~ N\in\N~~\forall ~n>N\qquad \frac{\log(\sum_{P\in\Pi_{n,\beta}}(\sharp P)^{
\frac{\log\beta}{\log\tau}})}{n\log\beta^{-1}}< \frac{\log(2\beta/\tau)}{\log(1/\tau)} .  \]
\end{proposition}
\proof
Fix a reciprocal of a PV number $\beta$. 
Consider the proof of Proposition 5.2 for the equal weighted Bernoulli measure $b$. Recall that we denote by
$x_{i}^{n}$ $i=1\dots\sharp(n)$ the distinct points of the
form $\sum_{k=0}^{n-1}\pm(1-\beta)\beta^{k}$ and by $m_{i}^{n}$ the
$b$ measure of corresponding element $P^{i}_{n}$ from the partition $\Pi_{n,\beta}$.\newline\newline
By the singularity of $b_{\beta}$ we have more than we used
in the proof of 5.2:\newline  
$\forall C\in(0,1)$ $\forall \epsilon >0~\exists $ disjoint intervals $(a_{1},b_{1}),\dots,(a_{u},b_{u})$ with 
\[ \sum_{l=1}^{u}(b_{l}-a_{l}) < \epsilon~~ \mbox{ and }~~ b_{\beta}(O)>C
~~\mbox{where}~~O:=\bigcup_{l=0}^{u}(a_{l},b_{l}).\]
By the same arguments we used in the proof of Proposition 5.2 we conclude:\newline
$\exists c>0$ $\forall C\in(0,1)$ $\forall \epsilon >0$ $\exists N=N(\epsilon,C)$ $\forall n\ge N$:
\[ \sum_{x_{i}^{n}\in \bar O}m_{i}^{n}>C \mbox{ and } \hat\sharp(n):=\mbox{card}\{x_{i}^{n}\in \bar O\}\le\epsilon c\beta^{-n}.\]
Since $m_{i}^{n}=b(P^{i}_{n})=\sharp P^{i}_{n}/2^{n}$, where  
$\sharp P$ denotes the number of
cylinder sets of length $n$ contained in $P$, it follows that
there is a subset $\hat\Pi_{n,\beta}$ of $\Pi_{n,\beta}$ with
$\hat\sharp(n)$ elements such that
\[ \sum_{P\in\hat\Pi_{n,\beta}}\sharp P \ge C 2^{n}\]
We estimate:
\[\sum_{P\in\Pi_{n,\beta}}(\sharp P)^{\log\beta/\log\tau}=
\sum_{P\in\hat\Pi_{n,\beta}} (\sharp P)^{\log\beta/\log\tau}
+\sum_{P\in\Pi_{n,\beta}\backslash \hat\Pi_{n,\beta} }(\sharp P)^{\log\beta/\log\tau}\]
\[
\le \hat\sharp(n)^{1-\log\beta/\log\tau}{(\sum_{P\in\hat\Pi_{n,\beta}}\sharp P)}^{\log\beta/\log\tau}\]\[+ (\sharp(n)-\hat\sharp(n))^{1-\log\beta/\log\tau}
{(\sum_{P\in\Pi_{n,\beta}\backslash \hat\Pi_{n,\beta} }\sharp P)}^{\log\beta/\log\tau}
\]
\[ \le (\epsilon c\beta^{-n})^{1-\log\beta/\log\tau} 2^{n \log\beta/\log\tau}+(c\beta^{-n})^{1-\log\beta/\log\tau}
((1-C)2)^{n\log\beta/\log\tau}\]
\[=\beta^{n(\log\beta/\log\tau-1)}2^{n \log\beta/\log\tau}
((\epsilon c)^{1-\log\beta/\log\tau} +c^{1-\log\beta/\log\tau}
(1-C)^{ \log\beta/\log\tau}).\]
Now choose $\epsilon$ and $C$ such that \[((\epsilon c)^{1-\log\beta/\log\tau} +c^{1-\log\beta/\log\tau}
(1-C)^{\log\beta/\log\tau})<1.\] 
For all $n\ge N(\epsilon,C)$ we have:
\[\frac{\log(\sum_{P\in\Pi_{n,\beta}}(\sharp P)^{
\frac{\log\beta}{\log\tau}})}{n\log\beta^{-1}}\] \[< \frac{\log(2\beta/\tau)}{\log(1/\tau)}+\frac{\log((\epsilon c)^{1-\log\beta/\log\tau} +c^{1-\log\beta/\log\tau}
(1-C)^{ \log\beta/\log\tau})}{n\log\beta^{-1}}.\]
The last term in this sum is negative and hence our proof is complete.\qed\newline\newline
{\bf Proof of 2.2}
From \cite{[PW]} we know that the box-counting dimension
of $\Lambda_{\beta,\tau}$ is given by $\log(2\beta/\tau)/\log(1/\tau)$. Thus Proposition 7.1 and 7.2  immediately imply $\dim_{H}\Lambda_{\beta,\tau}<\dim_{B}\Lambda_{\beta,\tau}$ if $\beta\in(0.5,1)$ is the reciprocal of a PV number. This is first statement of Theorem 2.2. Now the second statement remains to prove.
The following dimension formula for the Bernoulli measures $b^{p}_{\beta,\tau}$  
on $\Lambda_{\beta,\tau}$ is a 
corollary of Theorem II of \cite{[NE]}
\[ \dim_{H} b^{p}_{\beta,\tau}=\frac{p\log p+(1-p)\log(1-p)}{\log\tau}+(1-\frac{\log\beta}{\log\tau})\dim_{H}b^{p}_{\beta}. \]
Thus we have by Theorem 4.1 $\dim_{H}b^{p}_{\beta,\tau}<1$ for all $p\in(0,1)$ if $\beta\in(0.5,1)$ is the reciprocal of
a PV number and $\tau$ is small enough. But on the other hand we have $\dim_{H}\Lambda_{\beta,\tau}\ge 1$ since the projection of $\Lambda_{\beta,\tau}$ on the first coordinate axis is the whole interval $[-1,1]$. This proofs the second statement of our Theorem 2.2. \qed
\section*{Appendix A: General definitions and facts in dimension theory} 
We will here first define the most important quantities
in dimension theory and then collect some basic facts.
We refer to the book of Falconer \cite{[FA1]} and the book of
Pesin \cite{[PE1]} for a more detailed discussion of dimension theory.\newline
Let $q\in\N$ and $Z\subseteq\R^{q}$. For a real number $s>0$ we define the {\bf $s$-dimensional
Hausdorff measure} $H^{s}(Z)$ of $Z$ by 
\[ H^{s}(Z)=\lim_{\lambda \longrightarrow 0}\inf\{\sum_{i\in I}(\mbox{diam} U_{i})^{s}|Z\subseteq\bigcup_{i\in I}U_{i}\mbox{ and }\mbox{diam}(U_{i})\le\lambda\}\]
where $I$ is a countable index set.
The {\bf Hausdorff dimension} $\dim_{H}Z$ of $Z$ is given by 
\[ \dim_{H}Z=\sup\{s|H^{s}(Z)=\infty\}=\inf\{s|H^{s}(Z)=0\}.\]
Let $N_{\epsilon}(Z)$ be the minimal number of balls of radius
$\epsilon$ that are needed to cover $Z$. We define
the {\bf upper box-counting dimension} $\overline{\dim}_{B}$ resp. {\bf lower 
box-counting dimension} $\underline{\dim}_{B}$ of $Z$ by
\[ \overline{\dim}_{B}Z=\overline{\lim}_{\epsilon\longrightarrow 0} \frac{\log N_{\epsilon}(Z)}{-\log\epsilon}\qquad
\underline{\dim}_{B}Z=\underline{\lim}_{\epsilon\longrightarrow 0} \frac{\log N_{\epsilon}(Z)}{-\log\epsilon}.\]
If the limit it is called the {\bf box-counting dimension} $\dim_{B}$ of $Z$.
We remark that these quantities are not changed if we replace
$N_{\epsilon}(Z)$ by the minimal number of squares parallel to the axis with
side length $\epsilon$ that are needed to cover $Z$. Furthermore
we note that limit in the definition exists, if it exists
for some exponential decreasing sequence.\newline\newline
Now let $\mu$ be a Borel probability measure on
$\R^{q}$. We define the Hausdorff dimension of $\mu$ by
\[ \dim_{H}\mu=\inf\{\dim_{H}Z|\mu(Z)=1\}.\]
We introduce one more notion of dimension for the
measure $\mu$. 
Let $h_{\mu}(\epsilon)=\inf\{H_{\mu}(\Pi)|\Pi\mbox{ a partition with }\mbox{diam}\Pi\le\epsilon \}$ where $H_{\mu}(\Pi)$ is the usual
entropy of $\Pi$.
We define the {\bf upper R\'enyi dimension}  
$\overline{\dim}_{R}$ resp. {\bf lower  
R\'enyi dimension} $\underline{\dim}_{R}$ of $\mu$ by
\[ \overline{\dim}_{R}\mu=\overline{\lim}_{\epsilon\longrightarrow 0} \frac{ h_{\mu}(\epsilon)}{-\log\epsilon}\qquad
\underline{\dim}_{R}\mu=\underline{\lim}_{\epsilon\longrightarrow 0} \frac{h_{\mu}(\epsilon)}{-\log\epsilon}.\]
If the limit exists it is called  {\bf R\'enyi dimension} $\dim_{R}$ of $\mu$. 
The {\bf upper local dimension} $\overline{d}(x,\mu)$
resp. {\bf lower local dimension} $\underline{d}(x,\mu)$ of the measure
$\mu$ in a point $x$ is defined by 
\[ \overline{d}(x,\mu)=\overline{\lim}_{\epsilon\longrightarrow 0} \frac{\mu(B_{\epsilon}(x)) }{\log\epsilon}\qquad
\underline{d}(x,\mu)=\underline{\lim}_{\epsilon\longrightarrow 0} \frac{\mu(B_{\epsilon}(x))}{\log\epsilon}.\]
One basic fact we like to mention here is that dimensional theoretical quantities are not increased by projections or more general Lipschitz maps. This is immediate from the definitions. 
Basic relations between the dimensions
introduced here are stated in the following proposition.\newline\newline
{\bf Proposition A1} {\it For all $Z\subseteq \R^{q}$ and all Borel probability measures $\mu$ on $\R^{q}$ we have:\newline
(1) $\dim_{H}Z\le \underline{\dim}_{B}Z\le \overline{\dim}_{B}Z$ \newline
(2) $\underline{d}(x,\mu)\le c ~\mu-$almost everywhere $~ \Rightarrow~ 
\dim_{H}\mu\le c$.\newline
(3) $\underline{d}(x,\mu)\ge c ~\mu-$almost everywhere $~ \Rightarrow~
\dim_{H}\mu\ge c$ and ~$\underline{\dim}_{R}\mu\ge c$.\newline
(4) $\overline{d}(x,\mu)=\underline{d}(x,\mu)=c~\mu-$almost everywhere ~
$\Rightarrow~\dim_{H}\mu=\dim_{R}\mu=c$.}\newline\newline
The first inequality is obvious. A proof of the other statements is contained in the work of Young \cite{[YO]}. If the condition in part (4) holds,
the measure $\mu$ is called {\bf exact dimensional} and the common
value of the dimensions is denoted by $\dim\mu$. \newline\newline 
We need one other basic fact in our work which follows from Proposition 7.4 of \cite{[FA1]}..\newline\newline
{\bf Proposition A2}
{\it If $Z\subseteq\R^{q}$ and $I$ is an interval then
$\dim_{H}(Z\times I)=\dim_{H}+1$.}
\subsection*{Appendix B: Pisot-Vijayarghavan numbers}
A {\bf Pisot-Vijayarghavan number} (short: PV number) is
by definition the root of an algebraic equitation whose algebraic conjugates
lie all inside the unit circle in the complex plane.
Salem \cite{[SA]} showed that the set of PV numbers is a closed
subset of the reals and that $1$ is an isolated element.\newline
In our context we are interested in numbers $\beta \in(0.5,1)$
such that $\beta^{-1}$ is a PV number. We list
some examples including all reciprocals of PV numbers
with minimal polynomial of degree two and three and a sequence
of such numbers decreasing to $0.5$. 
\begin{center}
\begin{tabular}{|c|c|}
\hline
$x^{2}+x-1$ & $(\sqrt{5}-1)/2$\\        
\hline
$x^{3}+x^{2}+x-1$ & $0.5436898$\dots \\
\hline
$x^{3}+x^{2}-1$ & $0.754877$ \dots\\ 
\hline
$x^{3}+x-1$ & $0.6823278$\dots\\
\hline
$x^{3}-x^{2}+2x-1$ & $0.5698403$\dots \\
\hline
$x^{4}-x^{3}-1$ & $0.7244918$\dots \\
\hline
$x^{n}+x^{n-1}\dots+x-1$&~$r_{n}\longrightarrow 0.5$\\
\hline
\end{tabular}
\end{center}
{\bf Table 1:} Reciprocals of PV numbers \newline\newline
An important property of PV numbers is that their powers are near integers. More precise:\newline\newline
{\bf Proposition B1}{\it
If $\alpha$ is a PV number then there is a constant $0<\theta<1$
such that $||\alpha^{n}||_{\Z}\le  \theta^{n}$
$\forall n\ge0$ where $||.||_{\Z}$ denotes the distance to the
nearest integer.}\newline\newline
This statement can be found in \cite{[ER1]}. There is an another property of PV numbers that is of great importance for us.
For $\beta\in(0,1)$ we denote by $\sharp_{\beta}(n)$ the number
of distinct points of the for $\sum_{k=0}^{n-1}\pm\beta^{k}$ and by $\omega_{\beta}(n)$
the minimal distance between two of those points.\newline\newline 
{\bf Proposition B2}
{\it If $\beta\in(0.5,1)$ is the reciprocal of a PV number then
there are constants $\bar c>0$ and $\bar C>0$ such that $\omega_{\beta}(n)\ge \bar c\beta^{n}$ and
$\sharp_{\beta}(n)\ge \bar C \beta^{-n}$
holds for all
$n\ge 0$.}\newline\newline
For the first inequality we refer to Lemma 1.6 of \cite{[GA2]}. For the second inequality see formula (15) in \cite{[PU]}. 
Finally we like to mention that there is a whole
book about Pisot and Salem numbers \cite{[BDGPS]}. Certainly the reader
will find much more information about the role
of these numbers in algebraic number theory and Fourier analysis in this book than we provided here for our purposes. 


\begin{thebibliography}
\small
\bibitem{[AY]} J.C. Alexander, J.A. Yorke, Fat Baker s transformation's, Ergodic Thy. Dyn. Sys. 4, 1-23, 1984.
\bibitem{[BDGPS]} M.J. Bertin, A. Decomps-Guilloux, M. Grandet-Hugot, M. Pathiaux-Delefosse, J.P. Schreiber,
Pisot and Salem numbers, Birkhauser Verlag Basel, 1992.
\bibitem{[BPS]} L. Barreira, Ya. Pesin and J. Schmeling, Dimension and product structure of hyperbolic 
measures, Annals of Math., 149:3, 755-783, 1999.
\bibitem{[DGS]} M. Denker, C. Grillenberger, K.Sigmund, Ergodic Theory
on Compact Spaces, Lecture Notes in Math. 527, Springer Verlag Berlin, 1976.  
\bibitem{[ER1]} Erd\"os, On a family of symmetric Bernoulli convolutions,
, Amer. J. Math 61, 974-976, 1939.
\bibitem{[FA1]} K. Falconer, Fractal Geometry - Mathematical Foundations
and Applications, Wiley, New York, 1990.
\bibitem{[FA2]} K. Falconer,  The Hausdorff dimension of
self-affine fractals, Math. Proc. Camb. Phil. Soc. 103, 339-350, 1988.
\bibitem{[GA1]} A.M. Garsia, Entropy and singularity of infinite convolutions,
Pac. J. Math. 13, 1159-1169, 1963.
\bibitem{[GA2]} A.M. Garsia, Arithmetic properties of Bernoulli convolutions,
Trans. Amer. Math. Soc. 162, 409-432, 1962.
\bibitem{[HU]} J.E. Hutchinson, Fractals and self-similarity, Indiana Univ. Math. J. 30, 271-280, 1981.
\bibitem{[JW]} B. Jessen and A. Winter, Distribution functions and the Riemann zeta function,
Trans. Amer. Math. Soc. 38, 48-88, 1935.
\bibitem{[KH]} A. Katok and B. Hasselblatt, Introduction to Modern theory of
dynamical Systems, Cambridge University press, 1995.
\bibitem{[LY]} F. Ledrappier and L.-S. Young, The metric entropy of
diffeomorphism, Ann. Math. 122, 509-574, 1985.
\bibitem{[LP]} F. Ledrappier and A. Porzio, A dimension formula for Bernoulli convolutions, J. Stat. Phy. 76, no. 5/6, 1307-1326, 1994.
\bibitem{[MC]} C. McMullen, The Hausdorff dimension of
general Sierpinski carpets. Nagoya Math. J., 96, 1-9, 1984.
\bibitem{[MM]} A. Manning and H. McCluskey, Hausdorff dimension for horseshoes, Ergodic Thy. Dyn. Sys. 3, 251-260, 1983.
\bibitem{[NE]} J. Neunh\"auserer, Properties of some overlapping self-similar
and some self-affine measures, Schwerpunktprogramm
der deutschen Forschungsgemeinschaft: DANSE, Preprint 35/99; to appear
in: Acta Mathematica Hungarica 2002.
 \bibitem{[PE1]} Ya. Pesin, Dimension Theory in Dynamical Systems - 
Contemplary Views an Applications, University of Chicago Press, 1997. 
\bibitem{[PS2]} Y. Peres and B. Solomyak, Self-similar measures and intersection of Cantor sets, Trans. Amer. Math. Soc 350, no. 10, 4065-4087, 1998
\bibitem{[PSS]} Y. Peres, B. Solomyak and W. Schlag, Sixty years of Bernoulli
convolutions, Fractals and stochastic II, Progress in Probability 46, 95-106, Birkhauser, 2000.
\bibitem{[PU]} F. Przytychi and M. Urbanski, On Hausdorff dimension of some fractal sets, Studia. Math. 54, 218-228, 1989.
\bibitem{[PW]} M. Pollicott and H. Weiss, The dimension of self-affine limit
sets in the plane, J. Stat. Phys. 77, 841-860, 1994.
\bibitem{[SA]} R. Salem, A remarkable class of
algebraic integers, proof of a conjecture by Vijayarghavan,
Duke Math. J., 103-108, 1944.
\bibitem{[SCH]} J. Schmeling, A dimension formula for endomorphisms - The Belykh family,
Ergod. Th. Dyn. Sys 18, 1283-1309, 1998.
\bibitem{[SO1]} B. Solomyak, On the random series $\sum \pm \lambda^{i}$ (an Erd\"os problem), Ann. Math. 142, 1995.
\bibitem{[SO2]} B. Solomyak,  Measures and Dimensions for some Fractal Families, Proc. Cambridge Phil. Soc., 124/3, 531-546, 1998.
\bibitem{[SS]} B. Solomyak and K. Simon, Dimension of horseshoes in $\R^{3}$,
Ergodic Theory and Dyn. Sys. 19, 1345-1363, 1999.
\bibitem{[YO]} L.-S. Young, Dimension, entropy and Lyapunov exponents,
Ergod. Thy. Dyn. Sys. 2, 109-124, 1982.
\end{thebibliography}
\end{document}